\documentclass[twosided,reqno]{amsart}

\usepackage{amsmath}
\usepackage{amssymb}
\usepackage{amsthm}

\theoremstyle{theorem}
\newtheorem{theorem}{\scshape Theorem }[section]

\theoremstyle{definition}

\newtheorem{remark}{\scshape Remark}

\numberwithin{equation}{section}

\begin{document}

\title[Higher-order Cauchy of the first kind]{Higher-order Cauchy of the first kind and poly-Cauchy of the first kind mixed type polynomials}

\author{Dae San Kim}
\address{Department of Mathematics, Sogang University, Seoul 121-742, Republic of Korea.}
\email{dskim@sogang.ac.kr}

\author{Taekyun Kim}
\address{Department of Mathematics, Kwangwon University, Seoul 139-701, Republic of Korea}
\email{tkkim@kw.ac.kr}

\maketitle

\begin{abstract}
In this paper, we study higher-order Cauchy of the first kind and poly-Cauchy of the first kind mixed type polynomials with viewpoint of umbral calculus and give some interesting identities and formulae of those polynomials which are derived from umbral calculus.
\end{abstract}

\section{Introduction}

The polylogarithm factorial function is defined by

\begin{equation}\label{1}
Lif_{k}(t)=\sum_{n=0}^{\infty}\frac{t^{n}}{n!(n+1)^{k}},~~(k \in \mathbb{Z})~~(see~[13,14,15]).
\end{equation}

The poly-Cauchy polynomials of the first kind (of index $k$) are defined by the generating function to be

\begin{equation}\label{2}
\frac{Lif_{k}(log(1+t))}{(1+t)^{x}}=\sum_{n=0}^{\infty}C_{n}^{(k)}(x)\frac{t^{n}}{n!},~~(see~[13,14,15]).
\end{equation}
When $x=0$, $C_{n}^{(k)}=C_{n}^{(k)}(0)$ are called the  poly-Cauchy numbers.

In particular, for $k=1$, we note that

\begin{equation}\label{3}
Lif_{1}(log(1+t))=\frac{t}{\log (1+t)}=\sum_{n=0}^{\infty}C_{n}\frac{t^{n}}{n!},~~(see~[13,14]).
\end{equation}
where $C_{n}=C_{n}^{(1)}(0)$ are called the Cauchy numbers of the first kind.
The Cauchy numbers of the first kind with order $r$ are defined by the generating function to be

\begin{equation}\label{4}
\left(\frac{t}{\log(1+t)} \right)^{r}=\sum_{n=0}^{\infty}\mathbb{C}_{n}^{(r)}\frac{t^{n}}{n!},~~(see~[13,14,15]).
\end{equation}

Note that $\mathbb{C}_{n}^{(1)}=C_{n}$.
Let us consider higher-order Cauchy of the first kind and poly-Cauchy of the first kind mixed type polynomials as follows:

\begin{equation}\label{5}
\left(\frac{t}{\log(1+t)} \right)^{r}\frac{Lif_{k}(\log(1+t))}{(1+t)^{x}}=\sum_{n=0}^{\infty}A_{n}^{(r,k)}(x)\frac{t^{n}}{n!},
\end{equation}
where $r,k \in \mathbb{Z}$.
When $x=0$, $A_{n}^{(r,k)}=A_{n}^{(r,k)}(0)$ are called  higher-order Cauchy of the first kind and poly-Cauchy of the first kind mixed type numbers.

For $\lambda \neq 1 \in \mathbb{C}$, the Frobenius-Euler polynomials of order $\alpha$ are defined by the generating function to be

\begin{equation}\label{6}
\left(\frac{1-\lambda}{e^t-\lambda}\right)^{\alpha}e^{xt}=\sum_{n=0} ^{\infty}H_n ^{(\alpha)} (x|\lambda)\frac{t^n}{n!},~~{\text{ (see [9-12])}}.
\end{equation}

As is well known, the  Bernoulli polynomials of order $\alpha$ are given by

\begin{equation}\label{7}
\left(\frac{t}{e^t-1}\right)^{\alpha}e^{xt}=\sum_{n=0} ^{\infty}B_n ^{(\alpha)} (x)\frac{t^n}{n!},~~ {\text{ (see [1-8])}}.
\end{equation}
When $x=0$, $B_n ^{(\alpha)}=B_n ^{(\alpha)}(0)$ are called the  Bernoulli numbers of order $\alpha$.  The Stirling number of the first kind is defined by

\begin{equation}\label{8}
(x)_{n}=x(x-1)\cdots(x-n+1)=\sum_{l=0}^{n}S_{1}(n,l)x^{l},~~ {\text{ (see [16])}}.
\end{equation}

For $m \in \mathbb{Z}_{ \geq 0}$, the generating function of the Stirling number of the first kind is given by

\begin{equation}\label{9}
(\log (1+t))^{m}=m!\sum_{l=m} ^{\infty}S_{1}(l,m)\frac{t^l}{l!}=\sum_{l=0} ^{\infty}\frac{m!S_{1}(l+m,m)}{(l+m)!}t^{l+m}.
\end{equation}

Let $\mathbb{C}$ be the complex number field and let ${\mathcal{F}}$ be the set of all formal power series in the variable $t$:

\begin{equation}\label{10}
{\mathcal{F}}=\left\{ \left.f(t)=\sum_{k=0} ^{\infty} a_k\frac{ t^k}{k!}~\right|~ a_k \in {\mathbb{C}} \right\}.
\end{equation}

Let $\mathbb{P}=\mathbb{C}[x]$ and let $\mathbb{P}^{*}$ be the vector space of all linear functionals on $\mathbb{P}$.
$\left<\left. L~ ~\right|~ ~p(x)\right>$ is the action of the linear functional $L$ on the polynomial $p(x)$ with\\
 $\left< L+M~ \vert ~p(x)\right>=\left< L~ \vert ~p(x)\right>+\left< M~ \vert ~p(x)\right>$,
and $\left< cL~ \vert ~p(x)\right>=c \left< L~ \vert ~p(x)\right>$, where $c$ is a complex constant in $\mathbb{C}$.  For $f(t) \in \mathcal{F}$, let us define the linear functional on $\mathbb{P}$ by setting

\begin{equation}\label{11}
\left<f(t)|x^n \right>=a_k,~(n\geq 0), {\text{(see [1,7,16])}}.
\end{equation}

From (\ref{10}) and (\ref{11}), we note that

\begin{equation}\label{12}
\left<t^k | x^n \right>=n! \delta_{n,k},~(n,k \geq 0),
\end{equation}
where $\delta_{n,k}$ is the Kronecker symbol. (see [16, 17]). \\

Let $f_{L}(t)=\sum_{k=0}^{\infty}\frac{\left<L | x^{k} \right>}{k!}t^k$. Then by (\ref{12}), we get  $\left<f_{L}(t)|x^n \right>=\left<L | x^n \right>$.
So, $L=f_{L}(t)$.  The map $L \mapsto f_{L}(t)$ is a vector space isomorphism from $\mathbb{P}^{*}$ onto $\mathcal{F}$.  Henceforth, $\mathcal{F}$ denotes
both the algebra of formal power series in $t$ and the vector space of all linear functionals on $\mathbb{P}$, and so an element $f(t)$ of $\mathcal{F}$ will be thought of as both a formal power series and a linear functional.  We call $\mathcal{F}$ the umbral algebra and the umbral calculus is the study of umbral algebra. (see [16]).
The order $o(f(t))$ of a power series $f(t)(\neq 0)$ is the smallest integer $k$ for which the coefficient of $t^k$ does not vanish.  If $o(f(t))=1$, then $f(t)$ is called a delta series; if $o(f(t))=0$, then $g(t)$ is said to be an invertible seires.  For $f(t),g(t) \in \mathcal{F}$, let us assume that $f(t)$ is a delta series and $g(t)$ is an invertible series.
Then there exists a unique sequence $S_{n}(x)~(\deg S_{n}(x)=n)$ such that $\left<g(t)f(t)^k | S_n(x) \right>=n! \delta_{n,k}~for~n,k\geq 0$. The sequence $S_{n}(x)$ is called the Sheffer sequence for  $(g(t),f(t))$ which is
denoted by $S_n(x)\sim (g(t),f(t))$.  Let $f(t),g(t) \in \mathcal{F}$ and $p(x) \in \mathbb{P}$. Then we see that

\begin{equation}\label{13}
\left<f(t)g(t) | p(x) \right>=\left<f(t) | g(t)p(x) \right>= \left<g(t) | f(t)p(x) \right>,
\end{equation}
and

\begin{equation}\label{14}
f(t)=\sum_{k=0} ^{\infty} \left<f(t)|x^k\right>\frac{t^k}{k!},~p(x)=\sum_{k=0} ^{\infty}\left<t^k|p(x)\right> \frac{x^k}{k!},~(see[16]).
\end{equation}
From (\ref{14}), we note that

\begin{equation}\label{15}
t^{k}p(x)= p^{(k)}(x)= \frac{d^kp(x)}{dx^k},~and~e^{yt}p(x)=p(x+y).
\end{equation}

For $S_n(x) \sim (g(t), f(t))$, the generating function of $S_n(x)$ is given by

\begin{equation}\label{16}
\frac{1}{g(\bar{f}(t))}e^{x\bar{f}(t)}= \sum_{n=0}^{\infty}S_n(x)\frac{t^n}{n!},~for~all~x \in \mathbb{C},
\end{equation}
where $\bar{f}(t)$ is the compositional inverse of $f(t)$ with $\bar{f}(f(t))=t$.\\

From (\ref{5}), we observe that $A_{n}^{(r,k)}(x)$ is the Sheffer sequence for the pair\\
$\left( \left( \frac{te^{t}}{e^{t}-1} \right)^{r} \frac{1}{Lif_{k}(-t)}, e^{-t}-1 \right)$.
That is,

\begin{equation}\label{17}
A_{n}^{(r,k)}(x) \sim \left( \left( \frac{te^{t}}{e^{t}-1} \right)^{r} \frac{1}{Lif_{k}(-t)}, e^{-t}-1 \right).
\end{equation}
In [13], Komatsu considered the numbers $A_{n}^{(r,k)}$, which were denoted by $T_{r-1}^{(k)}(n)$.

Let $S_n(x) \sim (g(t), f(t))$.  Then we have

\begin{equation}\label{18}
f(t)S_n(x)=nS_{n-1}(x),~~(n \geq 0),~S_{n}(x)=\sum_{j=0}^{n}\frac{1}{j!} \langle \left. g(\bar{f}(t))^{-1}\bar{f}(t)^{j} ~\right|~ x^{n} \rangle x^{j},
\end{equation}

\begin{equation}\label{19}
S_{n}(x+y)=\sum_{j=0}^{n}\binom{n}{j}S_{j}(x)P_{n-j}(y),~where~p_{n}(x)=g(t)S_{n}(x),
\end{equation}
and

\begin{equation}\label{20}
S_{n+1}(x)=\left(x-\frac{g'(t)}{g(t)}\right)\frac{1}{f'(t)}S_n(x),~~(see~[16]).
\end{equation}

The transfer formula for $p_{n}(x) \sim (1,f(t))$, $q_{n}(x) \sim (1,g(t))$ is given by

\begin{equation}\label{21}
q_n(x)=x \left(\frac{f(t)}{g(t)}\right)^{n}x^{-1}p_{n}(x),~(n \geq 0)~(see~[16]).
\end{equation}

For $S_n(x) \sim (g(t),f(t))$,  $r_n(x) \sim (h(t),l(t))$, we have

\begin{equation}\label{22}
S_{n}(x)=\sum_{m=0}^{n}C_{n,m}r_m(x),~(n \geq 0),
\end{equation}

where

\begin{equation}\label{23}
C_{n,m}=\frac{1}{m!} \left< \left.\frac{h(\bar{f}(t))}{g(\bar{f}(t))}l(\bar{f}(t))^m ~\right|~ x^n \right>, ~~(see~[16]).
\end{equation}

In this paper, we study higher-order Cauchy of the first kind and poly-Cauchy of the first kind mixed type polynomials with viewpoint of umbral calculus. The purpose of this paper is to give some interesting identities and formulae of those polynomials which are derived from umbral calculus.

\section{Poly-Cauchy polynomials and Higher-order Cauchy polynomials}

By (\ref{17}), we see that

\begin{equation}\label{24}
 \left( \frac{te^{t}}{e^{t}-1} \right)^{r} \frac{1}{Lif_{k}(-t)} A_{n}^{(r,k)}(x) \sim (e^{-t}-1),
\end{equation}

and

\begin{equation}\label{25}
(-1)^{n}x^{(n)}=\sum_{m=0}^{n}(-1)^{m}S_{1}(n,m)x^{m} \sim (1,e^{-t}-1),
\end{equation}
where $x^{(n)}=x(x+1)\cdots(x+n-1)$.

Thus, from (\ref{24}) and (\ref{25}), we have

\begin{equation}\label{26}
 \left( \frac{te^{t}}{e^{t}-1} \right)^{r} \frac{1}{Lif_{k}(-t)} A_{n}^{(r,k)}(x)=(-1)^{n}x^{(n)} =\sum_{m=0}^{n}(-1)^{m}S_{1}(n,m)x^{m}.
\end{equation}

By (\ref{26}), we get

\begin{equation}\label{27}
\begin{split}
A_{n}^{(r,k)}(x)&=\left( \frac{e^{t}-1}{te^{t}} \right)^{r} Lif_{k}(-t) (-1)^{n}x^{(n)}\\
&=\sum_{m=0}^{n}(-1)^{m}S_{1}(n,m)\left( \frac{e^{t}-1}{te^{t}} \right)^{r} Lif_{k}(-t)x^{m}\\
&=\sum_{m=0}^{n}(-1)^{m}S_{1}(n,m)\sum_{l=0}^{m}\frac{(-1)^{l}(m)_{l}}{l!(l+1)^{k}} \left( \frac{e^{-t}-1}{-t} \right)^{r}x^{m-l}.
\end{split}
\end{equation}

It is well known that

\begin{equation}\label{28}
(e^t-1)^m=\sum_{l=0} ^{\infty}S_2(l+m,m)\frac{m!}{(l+m)!}t^{l+m},
\end{equation}
where $S_2(n,m)$ is the Stirling number of the second kind.

From (\ref{27}) and (\ref{28}), we have

\begin{equation}\label{29}
\begin{split}
A_{n}^{(r,k)}(x)&=\sum_{m=0}^{n}\sum_{l=0}^{m}\sum_{a=0}^{m-l}\frac{\binom{m}{l}\binom{m-l}{a}}{\binom{a+r}{r}(l+1)^{k}}S_1(n,m)S_2(a+r,r)(-x)^{m-l-a}\\
&=\sum_{m=0}^{n}\left\{ \sum_{l=0}^{m}\sum_{j=0}^{m-l}\frac{\binom{m}{l}\binom{m-l}{j}}{\binom{m-l-j+r}{r}(l+1)^{k}}S_1(n,m)S_2(m-l-j+r,r)\right\}(-x)^{j}\\
&=\sum_{j=0}^{n}\left\{ \sum_{m=j}^{n}\sum_{l=0}^{m-j}\frac{\binom{m}{l}\binom{m-l}{j}}{\binom{m-l-j+r}{r}(l+1)^{k}}S_1(n,m)S_2(m-l-j+r,r)\right\}(-x)^{j}\\
\end{split}
\end{equation}
where $r \in \mathbb{Z}_{\geq 0}=\mathbb{N}\cup \{0\}$.

Therefore, by (\ref{29}), we obtain the following theorem.

\begin{theorem}\label{thm1}
For $n,r \geq 0$, we have
\begin{equation*}
A_{n}^{(r,k)}(x)=\sum_{0 \leq j \leq n }\left\{ \sum_{m=j}^{n}\sum_{l=0}^{m-j}\frac{\binom{m}{l}\binom{m-l}{j}}{\binom{m-l-j+r}{r}(l+1)^{k}}S_1(n,m)S_2(m-l-j+r,r)\right\}(-x)^{j}.
\end{equation*}
\end{theorem}

From (\ref{17}) and (\ref{18}), we have

\begin{equation}\label{30}
\begin{split}
&A_{n}^{(r,k)}(x)\\
&=\sum_{j=0}^{n} \frac{1}{j!} \left< \left.\left( \frac{t}{\log (1+t)} \right)^{r} Lif_{k}(\log (1+t)) (-\log (1+t))^{j} ~\right|~ x^{n} \right> x^{j}\\
&=\sum_{j=0}^{n}\sum_{l=0}^{n-j}(-1)^{j}\binom{n}{l+j}S_1(l+j,j) \left<\left. Lif_{k}(\log (1+t))\left( \frac{t}{\log (1+t)} \right)^{r}  ~\right|~ x^{n-l-j} \right> x^{j}\\
&=\sum_{j=0}^{n}\sum_{l=0}^{n-j}(-1)^{j}\binom{n}{l+j}S_1(l+j,j)\sum_{a=0}^{n-l-j}B_{a}^{(a-r+1)}(1)\binom{n-l-j}{a} \\
&\qquad\qquad\qquad\qquad\qquad\qquad\qquad \times \left<\left. Lif_{k}(\log (1+t)) ~\right|~ x^{n-l-j-a} \right> x^{j}\\
&=\sum_{j=0}^{n}\left\{\sum_{l=0}^{n-j}\sum_{a=0}^{n-l-j}(-1)^{j}\binom{n}{l+j}\binom{n-l-j}{a}S_1(l+j,j)B_{a}^{(a-r+1)}(1)C_{n-j-l-a}^{(k)}\right\}x^{j}.\\
\end{split}
\end{equation}

Therefore, by (\ref{30}), we obtain the following theorem.

\begin{theorem}\label{thm2}
For $r,k \in \mathbb{Z}$, and $n \in \mathbb{Z}_{\geq 0}$, we have
\begin{equation*}
\begin{split}
&A_{n}^{(r,k)}(x)\\
&=\sum_{j=0}^{n}\left\{\sum_{l=0}^{n-j}\sum_{a=0}^{n-l-j}(-1)^{j}\binom{n}{l+j}\binom{n-l-j}{a}S_1(l+j,j)B_{a}^{(a-r+1)}(1)C_{n-j-l-a}^{(k)}\right\}x^{j}.\\
\end{split}
\end{equation*}
\end{theorem}

As is knwon, the Narumi polynomials of order $r$ are given by $N_{n}^{(r)}(x) \sim $\\
$\left( \left(\frac{e^{t}-1}{t} \right)^{r}, (e^{t}-1) \right)$.
Thus, we note that

\begin{equation}\label{31}
\left( \frac{t}{\log (1+t)} \right)^{-r}(1+t)^{x}= \sum_{n=0}^{\infty}N_{n}^{(r)}(x)\frac{t^{n}}{n!}.
\end{equation}

Indeed, we see that $N_{n}^{(r)}(x)=B_{n}^{(n+r+1)}(x+1)$.

From (\ref{30}) and (\ref{31}), we can derive the following equation:

\begin{equation}\label{32}
A_{n}^{(r,k)}(x)=\sum_{j=0}^{n}\left\{\sum_{l=0}^{n-j}\sum_{a=0}^{n-l-j}(-1)^{j}\binom{n}{l+j}\binom{n-l-j}{a}S_1(l+j,j)N_{a}^{(-r)}C_{n-j-l-a}^{(k)}\right\}x^{j},
\end{equation}
where $N_{a}^{(r)}=N_{a}^{(r)}(0)$ are called the Narumi numbers of order $r$.

The Bernoulli polynomials of the second kind are defined by the generating function to be

\begin{equation}\label{33}
\frac{t}{\log (1+t)}(1+t)^{x}= \sum_{n=0}^{\infty}b_{n}(x)\frac{t^{n}}{n!},~~(see~[16]).
\end{equation}

From (\ref{30}) and (\ref{33}), we note that

\begin{equation}\label{34}
\begin{split}
A_{n}^{(r,k)}(x)&=\sum_{j=0}^{n}\left\{\sum_{l=0}^{n-j}\sum_{a=0}^{n-l-j}\sum_{a_{1}+\cdots +a_{r}=a}(-1)^{j} \binom{n}{l+j}\binom{n-l-j}{a} \binom{a}{a_{1},\cdots ,a_{r}}\right.\\
&\left. \times S_{1}(l+j,j)\left( \Pi_{i=1}^{r}b_{a_{i}}\right)C_{n-l-j-a}^{(k)}\right\}x^{j}.\\
\end{split}
\end{equation}

From (\ref{17}),(\ref{19}) and (\ref{26}), we note that

\begin{equation}\label{35}
A_{n}^{(r,k)}(x+y)=\sum_{j=0}^{n}(-1)^{n-j}\binom{n}{j}A_{j}^{(r,k)}(x)y^{(n-j)},
\end{equation}

and, by (\ref{15}) and (\ref{18}), we get

\begin{equation}\label{36}
nA_{n-1}^{(r,k)}(x)=(e^{-t}-1)A_{n}^{(r,k)}(x)=A_{n}^{(r,k)}(x-1)-A_{n}^{(r,k)}(x).
\end{equation}

By (\ref{17}) and (\ref{20}), we get

\begin{equation}\label{37}
\begin{split}
A_{n+1}^{(r,k)}(x)&=\left( \frac{g'(t)}{g(t)}-x \right)e^{t}A_{n}^{(r,k)}(x)\\
&=e^{t}\frac{g'(t)}{g(t)}A_{n}^{(r,k)}(x)-xA_{n}^{(r,k)}(x+1)\\
&=r \frac{e^{t}-1-t}{t^2}\frac{te^{t}}{e^{t}-1}A_{n}^{(r,k)}(x)+e^{t}\frac{Lif_{k}'(-t)}{Lif_{k}(-t)}A_{n}^{(r,k)}(x)-xA_{n}^{(r,k)}(x+1).\\
\end{split}
\end{equation}

From (\ref{26}), we note that

\begin{equation}\label{38}
A_{n}^{(r,k)}(x)=\sum_{m=0}^{n}(-1)^{m}S_{1}(n,m)\left( \frac{e^{t}-1}{te^{t}} \right)^{r}Lif_{k}(-t) x^{m},
\end{equation}

\begin{equation}\label{39}
\frac{1}{Lif_{k}(-t)}A_{n}^{(r,k)}(x)=\sum_{m=0}^{n}(-1)^{m}S_{1}(n,m)\left( \frac{e^{t}-1}{te^{t}} \right)^{r} x^{m}.
\end{equation}

By (\ref{38}), we get

\begin{equation}\label{40}
\begin{split}
&r\left( \frac{e^{t}-1-t}{t^{2}} \right)\left( \frac{te^{t}}{e^{t}-1} \right)A_{n}^{(r,k)}(x)\\
&=r\sum_{m=0}^{n}(-1)^{m}S_{1}(n,m)\frac{e^{t}-1-t}{t^{2}}\left( \frac{te^{t}}{e^{t}-1} \right)^{1-r}Lif_{k}(-t) x^{m}\\
&=r\sum_{m=0}^{n}(-1)^{m}S_{1}(n,m)\sum_{l=0}^{m}\frac{(-1)^{l}(m)_{l}}{l!(l+1)^{k}}\left( \frac{te^{t}}{e^{t}-1} \right)^{1-r}\sum_{a=0}^{m-l}\frac{t^{a}x^{m-l}}{(a+2)!}\\
&=r\sum_{m=0}^{n}(-1)^{m}S_{1}(n,m)\sum_{l=0}^{m}\frac{(-1)^{l}(m)_{l}}{l!(l+1)^{k}}\sum_{a=0}^{m-l}\frac{(m-l)_{a}}{(a+2)!}\left( \frac{-t}{e^{-t}-1} \right)^{1-r}x^{m-l-a}\\
&=r\sum_{m=0}^{n}\sum_{l=0}^{m}\sum_{a=0}^{m-l}\frac{(-1)^{a}\binom{m}{l}\binom{m-l}{a}}{(a+2)(a+1)(l+1)^{k}}S_{1}(n,m)B_{m-l-a}^{(1-r)}(-x),\\
\end{split}
\end{equation}

and

\begin{equation}\label{41}
\begin{split}
&e^{t}\frac{Lif_{k}'(-t)}{Lif_{k}(-t)}A_{n}^{(r,k)}(x)=e^{t}Lif_{k}'(-t)\left(\frac{1}{Lif_{k}(-t)}A_{n}^{(r,k)}(x)\right)\\
&=e^{t}Lif_{k}'(-t)\sum_{m=0}^{n}(-1)^{m}S_{1}(n,m)\left( \frac{e^{t}-1}{te^{t}} \right)^{r} x^{m}\\
&=\sum_{m=0}^{n}(-1)^{m}S_{1}(n,m)e^{t}\left( \frac{e^{t}-1}{te^{t}} \right)^{r}\sum_{a=0}^{m}\frac{(-1)^{a}}{a!(a+2)^{k}}t^{a}x^{m}\\
&=\sum_{m=0}^{n}(-1)^{m}S_{1}(n,m)\sum_{a=0}^{m}\frac{(-1)^{a}}{a!(a+2)^{k}}(m)_{a}e^{t}\left( \frac{-t}{e^{-t}-1} \right)^{-r}x^{m-a}\\
&=\sum_{m=0}^{n}\sum_{a=0}^{m}\binom{m}{a}\frac{S_{1}(n,m)}{(a+2)^{k}}B_{m-a}^{(-r)}(-x-1).\\
\end{split}
\end{equation}

Therefore, by (\ref{37}), (\ref{40}) and (\ref{41}), we obtain the following theorem.

\begin{theorem}\label{thm3}
For $r,k \in \mathbb{Z}$, and $n \geq 0$, we have
\begin{equation*}
\begin{split}
A_{n+1}^{(r,k)}(x)&=-xA_{n}^{(r,k)}(x+1)+r\sum_{m=0}^{n}\sum_{l=0}^{m}\sum_{a=0}^{m-l}\frac{(-1)^{a}\binom{m}{l}\binom{m-l}{a}}{(a+2)(a+1)(l+1)^{k}}\\
&\times S_{1}(n,m)B_{m-l-a}^{(1-r)}(-x)+\sum_{m=0}^{n}\sum_{a=0}^{m}\binom{m}{a}\frac{S_{1}(n,m)}{(a+2)^{k}}B_{m-a}^{(-r)}(-x-1).\\
\end{split}
\end{equation*}
\end{theorem}

By (\ref{12}), we easily see that
\begin{equation}\label{42}
\begin{split}
&A_{n}^{(r,k)}(y)= \left<\left. \sum_{l=0}^{\infty}A_{l}^{(r,k)}(y)\frac{t^{l}}{l!} ~\right|~ x^{n} \right>\\
&=\left< \left.\left( \frac{t}{\log (1+t)} \right)^{r} Lif_{k}(\log(1+t))(1+t)^{-y} ~\right|~ x^{n} \right>\\
&=\left< \left.\left( \frac{t}{\log (1+t)} \right)^{r} Lif_{k}(\log(1+t))(1+t)^{-y} ~\right|~ xx^{n-1} \right>\\
&=\left<\left. \partial_{t} \left\{ \left( \frac{t}{\log (1+t)} \right)^{r} Lif_{k}(\log(1+t))(1+t)^{-y}\right\} ~\right|~ x^{n-1} \right>\\
&=\left<\left. \partial_{t} \left( \left( \frac{t}{\log (1+t)} \right)^{r} \right)Lif_{k}(\log(1+t))(1+t)^{-y} ~\right|~ x^{n-1} \right>\\
& \qquad + \left< \left. \left( \frac{t}{\log (1+t)} \right)^{r} (\partial_{t}Lif_{k}(\log(1+t)))(1+t)^{-y} ~\right|~ x^{n-1} \right>\\
& \qquad + \left< \left. \left( \frac{t}{\log (1+t)} \right)^{r}Lif_{k}(\log(1+t))(\partial_{t}(1+t)^{-y}) ~\right|~ x^{n-1} \right>\\
&=-yA_{n-1}^{(r,k)}(y+1)+\left<\left. \left(\partial_{t}\left( \frac{t}{\log (1+t)} \right)^{r}\right)Lif_{k}(\log(1+t))(1+t)^{-y} ~\right|~ x^{n-1} \right>\\
&\qquad +\left< \left.\left( \frac{t}{\log (1+t)} \right)^{r}(\partial_{t}Lif_{k}(\log(1+t)))(1+t)^{-y}~\right|~ x^{n-1} \right>.\\
\end{split}
\end{equation}

Now, we observe that

\begin{equation}\label{43}
\begin{split}
&\left<\left. \left(\partial_{t}\left( \frac{t}{\log (1+t)} \right)^{r}\right)Lif_{k}(\log(1+t))(1+t)^{-y} ~\right|~ x^{n-1} \right>\\
&=r\sum_{l=0}^{n-1}\sum_{a=0}^{l}(-1)^{n-a}\frac{(n-1-l)!(l-a)!}{l-a+2}\binom{n-1}{l}\binom{l}{a}A_{a}^{(r+1,k)}(y)\\
&+r\sum_{l=0}^{n-1}(-1)^{n-1-l}(n-1-l)!\binom{n-1}{l}A_{l}^{(r,k)}(y),\\
\end{split}
\end{equation}

and
\begin{equation}\label{44}
\begin{split}
&\left<\left(\frac{t}{\log(1+t)}\right)^{r}(\partial _{t}Lif_{k}(\log (1+t))(1+t)^{-y}|x^{n-1} \right>\\
&=\frac{1}{n}(A_{n}^{(r+1,k-1)}(y+1)-A_{n}^{(r+1,k)}(y+1)).\\
\end{split}
\end{equation}

Therefore, by (\ref{42}),(\ref{43}) and (\ref{44}), we obtain the following theorem.

\begin{theorem}\label{thm4}
For $r,k \in \mathbb{Z}$ and $n \geq 0$, we have
\begin{equation*}
\begin{split}
A_{n}^{(r,k)}(x)&=-xA_{n-1}^{(r,k)}(x+1)+r\sum_{l=0}^{n-1}\sum_{a=0}^{l}(-1)^{n-a}\frac{(n-1-l)!(l-a)!}{l-a+2}\binom{n-1}{l}\\
&\times \binom{l}{a}A_{n}^{(r+1,k)}(x)+r\sum_{l=0}^{n-1}(-1)^{n-l-1}(n-l-1)!\binom{n-1}{l}A_{l}^{(r,k)}(x)\\
&+\frac{1}{n}(A_{n}^{(r+1,k-1)}(x+1)-A_{n}^{(r+1,k)}(x+1)).\\
\end{split}
\end{equation*}
\end{theorem}
Here we compute $\left<\left(\frac{t}{\log(1+t)}\right)^{r}Lif_{k}(\log (1+t))(\log (1+t))^{m}|x^{n} \right>$ in two different ways.

On the one hand, we have

\begin{equation}\label{45}
\begin{split}
&\left<\left(\frac{t}{\log(1+t)}\right)^{r}Lif_{k}(\log (1+t))(\log (1+t))^{m}|x^{n} \right>\\
&=\sum_{l=0}^{n-m}\frac{m!}{(l+m)!}S_{1}(l+m,m)(n)_{l+m}\left<\left(\frac{t}{\log(1+t)}\right)^{r}Lif_{k}(\log (1+t))|x^{n-l-m} \right>\\
&=\sum_{l=0}^{n-m}m!\binom{n}{l+m}S_{1}(l+m,m)A_{n-l-m}^{(r,k)}\\
&=\sum_{l=0}^{n-m}m!\binom{n}{l}S_{1}(n-l,m)A_{l}^{(r,k)}.\\
\end{split}
\end{equation}

On the other hand, we get

\begin{equation}\label{46}
\begin{split}
&\left<\left(\frac{t}{\log(1+t)}\right)^{r}Lif_{k}(\log (1+t))(\log (1+t))^{m}|x^{n} \right>\\
&=\left<\left(\frac{t}{\log(1+t)}\right)^{r}Lif_{k}(\log (1+t))(\log (1+t))^{m}|xx^{n-1} \right>\\
&=\left<\partial_{t}\left\{\left(\frac{t}{\log(1+t)}\right)^{r}Lif_{k}(\log (1+t))(\log (1+t))^{m}\right\}|x^{n-1} \right>\\
&=\left<\left( \partial_{t}\left(\frac{t}{\log(1+t)}\right)^{r}\right)Lif_{k}(\log (1+t))(\log (1+t))^{m}|x^{n-1} \right>\\
&\qquad +\left<\left(\frac{t}{\log(1+t)}\right)^{r}\left(\partial_{t}Lif_{k}(\log (1+t))\right)(\log (1+t))^{m}|x^{n-1} \right>\\
&\qquad +\left<\left(\frac{t}{\log(1+t)}\right)^{r}Lif_{k}(\log (1+t))\left(\partial_{t}(\log (1+t))^{m}\right)|x^{n-1} \right>.\\
\end{split}
\end{equation}

Now, we observe that

\begin{equation}\label{47}
\begin{split}
&\left<\left( \partial_{t}\left(\frac{t}{\log(1+t)}\right)^{r}\right)Lif_{k}(\log (1+t))(\log (1+t))^{m}|x^{n-1} \right>\\
&=r\sum_{l=0}^{n-1-m}m!\binom{n-1}{l}S_{1}(n-l-1,m)A_{l}^{(r,k)}(1)\\
&+r\sum_{l=0}^{n-1-m}\sum_{a=0}^{l}(-1)^{l-a+1}\frac{m!(l-a)!}{l-a+2}\binom{n-1}{l}\binom{l}{a}S_{1}(n-1-l,m)A_{a}^{(r+1,k)}(1),\\
\end{split}
\end{equation}

\begin{equation}\label{48}
\begin{split}
&\left<\left(\frac{t}{\log(1+t)}\right)^{r}\left( \partial_{t}Lif_{k}(\log (1+t))(\log (1+t))^{m}\right)|x^{n-1} \right>\\
&=\sum_{l=0}^{n-m}(m-1)!\binom{n-1}{l}S_{1}(n-l-1,m-1)\{A_{l}^{(r,k-1)}(1)-A_{l}^{(r,k)}(1)\},\\
\end{split}
\end{equation}

and

\begin{equation}\label{49}
\begin{split}
&\left<\left(\frac{t}{\log(1+t)}\right)^{r}Lif_{k}(\log (1+t))\left( \partial_{t}(\log (1+t))^{m}\right)|x^{n-1} \right>\\
&=m!\sum_{l=0}^{n-m}\binom{n-1}{l}S_{1}(n-l-1,m-1)A_{l}^{(r,k)}(1),\\
\end{split}
\end{equation}
where $n-1 \geq m \geq 1$.

From (\ref{46}),(\ref{47}),(\ref{48}) and (\ref{49}), we have

\begin{equation}\label{50}
\begin{split}
&\left<\left(\frac{t}{\log(1+t)}\right)^{r}Lif_{k}(\log (1+t))(\log (1+t))^{m}|x^{n} \right>\\
&=r\sum_{l=0}^{n-1-m}\sum_{a=0}^{l}(-1)^{l-a+1}\frac{m!(l-a)!}{l-a+2}\binom{n-1}{l}\binom{l}{a}S_{1}(n-1-l,m)A_{a}^{(r+1,k)}(1)\\
&+r\sum_{l=0}^{n-1-m}m!\binom{n-1}{l}S_{1}(n-l-1,m)A_{l}^{(r,k)}(1)+\sum_{l=0}^{n-m}(m-1)!\binom{n-1}{l}\\
& \times S_{1}(n-l-1,m-1)A_{l}^{(r,k-1)}(1)-\sum_{l=0}^{n-m}(m-1)!\binom{m-1}{l}S_{1}(n-l-1,m-1)\\
&\times A_{l}^{(r,k)}(1)+m!\sum_{l=0}^{n-m}\binom{n-1}{l}S_{1}(n-l-1,m-1)A_{l}^{(r,k)}(1).\\
\end{split}
\end{equation}

Therefore, by (\ref{45}) and (\ref{50}), we obtain the following theorem.

\begin{theorem}\label{thm5}
For $n-1 \geq m \geq 1$, we have
\begin{equation*}
\begin{split}
&\sum_{l=0}^{n-m}\binom{n}{l}S_{1}(n-l,m)A_{l}^{(r,k)}\\
&=r\sum_{l=0}^{n-1-m}\sum_{a=0}^{l}(-1)^{l-a+1}\frac{(l-a)!}{(l-a+2)}\binom{n-1}{l}\binom{l}{a}S_{1}(n-1-l,m)A_{a}^{(r+1,k)}(1)\\
&+r\sum_{l=0}^{n-1-m}\binom{n-1}{l}S_{1}(n-l-1,m)A_{l}^{(r,k)}(1)+\frac{1}{m}\sum_{l=0}^{n-m}\binom{n-1}{l}\\
&\times S_{1}(n-l-1,m-1)A_{l}^{(r,k)}(1)+\left(1-\frac{1}{m}\right)\sum_{l=0}^{n-m}\binom{n-1}{l}S_{1}(n-l-1,m-1)\\
& \times A_{l}^{(r,k)}(1).\\
\end{split}
\end{equation*}
\end{theorem}

\begin{remark}\label{rmk}
It is known that
\begin{equation}\label{51}
\frac{d}{dx}S_{n}(x)=\sum_{l=0}^{n-1}\binom{n}{l} \left<\left. \bar{f}(t) ~\right|~ x^{n-l} \right> S_{l}(x),~~(see[16]),
\end{equation}
\end{remark}
where $S_n(x)\sim (g(t),f(t))$.

From (\ref{17}) and (\ref{51}), we have

\begin{equation}\label{52}
\begin{split}
\frac{d}{dx}A_{n}^{(r,k)}(x)&=\sum_{l=0}^{n-1}\binom{n}{l} \left< \left.-\log (1+t) ~\right|~ x^{n-l} \right> A_{l}^{(r,k)}(x)\\
&=-\sum_{l=0}^{n-1}\binom{n}{l} \left< \left.\sum_{m=0}^{\infty}\frac{(-1)^{m}}{m+1}t^{m+1} ~\right|~ x^{n-l} \right> A_{l}^{(r,k)}(x)\\
&=-\sum_{l=0}^{n-1}\binom{n}{l}\sum_{m=0}^{\infty}\frac{(-1)^{m}}{m+1} \left< \left. t^{m+1} ~\right|~ x^{n-l} \right> A_{l}^{(r,k)}(x)\\
&=(-1)^{n+1}n!\sum_{l=0}^{n-1}\frac{(-1)^{l+1}}{(n-l)l!}A_{l}^{(r,k)}(x).\\
\end{split}
\end{equation}

For $A_{n}^{(r,k)}(x) \sim \left(\left(\frac{te^{t}}{e^{t}-1}\right)^{r}\frac{1}{Lif_{k}(-t)},e^{-t}-1\right) $ and $B_{n}^{(s)}(x) \sim \left(\left(\frac{e^{t}-1}{t}\right)^{s},t \right)$,\\
 $(s \geq 0)$, let us assume that

\begin{equation}\label{53}
A_{n}^{(r,k)}(x)=\sum_{m=0}^{n}C_{n,m}B_{m}^{(s)}(x).
\end{equation}

Then, by (\ref{23}), we get

\begin{equation}\label{54}
\begin{split}
&C_{n,m}=\frac{1}{m!}\left<\left(\frac{t}{(1+t)\log(1+t)}\right)^{r+s}(1+t)^{r}Lif_{k}(\log (1+t))(-\log (1+t))^{m}|x^{n} \right>\\
&=\frac{(-1)^{m}}{m!}\left<\left(\frac{t}{\log(1+t)}\right)^{r+s}Lif_{k}(\log (1+t))(1+t)^{-s}|(\log (1+t))^{m}x^{n} \right>\\
&=(-1)^{m}\sum_{l=0}^{n-m}\binom{n}{l+m}S_{1}(l+m,m)\\
&\quad\times \left<\left(\frac{t}{\log(1+t)}\right)^{r+s}Lif_{k}(\log (1+t))(1+t)^{-s}|x^{n-l-m} \right>\\
&=(-1)^{m}\sum_{l=0}^{n-m}\binom{n}{l+m}S_{1}(l+m,m)A_{n-l-m}^{(r+s,k)}(s)\\
&=(-1)^{m}\sum_{l=0}^{n-m}\binom{n}{l}S_{1}(n-l,m)A_{l}^{(r+s,k)}(s).\\
\end{split}
\end{equation}

Therefore, by (\ref{53}) and (\ref{54}), we obtain the following theorem.

\begin{theorem}\label{thm6}
For $n,s \geq 0$, we have
\begin{equation*}
\begin{split}
A_{n}^{(r,k)}(x)=\sum_{m=0}^{n}\left\{(-1)^{m}\sum_{l=0}^{n-m}\binom{n}{l}S_{1}(n-l,m)A_{l}^{(r+s,k)}(s)\right\}B_{m}^{(s)}(x).
\end{split}
\end{equation*}
\end{theorem}

Let us consider the following two Sheffer sequences:

\begin{equation}\label{55}
A_{n}^{(r,k)}(x) \sim \left(\left(\frac{te^{t}}{e^{t}-1}\right)^{r}\frac{1}{Lif_{k}(-t)},e^{-t}-1\right),
\end{equation}

and

\begin{equation}\label{56}
H_{n}^{(s)}(x|\lambda) \sim \left(\left(\frac{e^{t}-\lambda}{1-\lambda}\right)^{s},t \right),~~(s \geq 0).
\end{equation}

Suppose that

\begin{equation}\label{57}
A_{n}^{(r,k)}(x)=\sum_{m=0}^{n}C_{n,m}H_{m}^{(s)}(x|\lambda)~~(s \geq 0).
\end{equation}

By (\ref{23}), we get

\begin{equation}\label{58}
\begin{split}
&C_{n,m}\\
&=\frac{(-1)^{m}}{m!(1-\lambda)^{s}}\\
&\quad \times \left<\left(\frac{t}{\log(1+t)}\right)^{r}Lif_{k}(\log (1+t))(1+t)^{-s}(1-\lambda(1+t))^{s}|(\log (1+t))^{m}x^{n} \right>\\
&=\frac{(-1)^{m}}{(1-\lambda)^{s}}\sum_{l=0}^{n-m}\sum_{a=0}^{s}(-\lambda)^{a}\binom{n}{l+m}\binom{s}{a}S_{1}(l+m,m)A_{n-l-m}^{(r,k)}(s-a)\\
&=\frac{(-1)^{m}}{(1-\lambda)^{s}}\sum_{l=0}^{n-m}\sum_{a=0}^{s}(-\lambda)^{a}\binom{n}{l}\binom{s}{a}S_{1}(n-l,m)A_{l}^{(r,k)}(s-a).\\
\end{split}
\end{equation}

Therefore, by (\ref{57}) and (\ref{58}), we obtain the following theorem.

\begin{theorem}\label{thm7}
For $n,s \geq 0$, $r,k \in \mathbb{Z}$, we have
\begin{equation*}
\begin{split}
&A_{n}^{(r,k)}(x)=\frac{1}{(1-\lambda)^{s}}\sum_{m=0}^{n}\left\{(-1)^{m}\sum_{l=0}^{n-m}\sum_{a=0}^{s}(-\lambda)^{a}\binom{n}{l}\binom{s}{a}S_{1}(n-l,m)A_{l}^{(r,k)}(s-a)
\right\}\\
&\qquad \qquad\quad \times H_{m}^{(s)}(x|\lambda).\\
\end{split}
\end{equation*}
\end{theorem}

Finally, we consider

\begin{equation}\label{59}
A_{n}^{(r,k)}(x) \sim \left(\left(\frac{te^{t}}{e^{t}-1}\right)^{r}\frac{1}{Lif_{k}(-t)},e^{-t}-1\right),
\end{equation}

and

\begin{equation}\label{60}
x^{(n)} \sim (1,e^{-t}-1).
\end{equation}

Let us assume that

\begin{equation}\label{61}
A_{n}^{(r,k)}(x)=\sum_{m=0}^{n}C_{n,m}x^{(m)}.
\end{equation}

Then, by (\ref{23}), we get

\begin{equation}\label{62}
\begin{split}
C_{n,m}&=\frac{(-1)^{m}}{m!}\left<\left(\frac{t}{\log(1+t)}\right)^{r}Lif_{k}(\log (1+t))|t^{m}x^{n} \right>\\
&=(-1)^{m}\frac{(n)_{m}}{m!}\left<\left(\frac{t}{\log(1+t)}\right)^{r}Lif_{k}(\log (1+t))|x^{n-m} \right>\\
&=(-1)^{m}\binom{n}{m}A_{n-m}^{(r,k)}.\\
\end{split}
\end{equation}

Therefore, by (\ref{61}) and (\ref{62}), we obtain the following theorem.

\begin{theorem}\label{thm8}
For $n\geq 0$, $r,k \in \mathbb{Z}$, we have
\begin{equation*}
\begin{split}
A_{n}^{(r,k)}(x)=\sum_{m=0}^{n}(-1)^{m}\binom{n}{m}A_{n-m}^{(r,k)}x^{(m)},
\end{split}
\end{equation*}
\end{theorem}
where $x^{(m)}=x(x+1) \cdots (x+m-1).$


\bigskip
ACKNOWLEDGEMENTS. This work was supported by the National Research Foundation of Korea(NRF) grant funded by the Korea government(MOE)\\
(No.2012R1A1A2003786 ).
\bigskip



\begin{thebibliography}{99}


\bibitem {01} S. Araci and M. Acikgoz, {\it A note on the Frobenius-Euler numbers and polynomials associated with Bernstein polynomials}, Adv. Stud. Contemp. Math., ${\mathbf{22}}$  (2012),  no.3,  399-406.

\bibitem {02}A. Bayad, T. Kim, {\it Identities involving values of Bernstein, $q$-Bernoulli, and $q$-Euler polynomials}, Russ, J. Math. Phys. 18 (2011), no. 2, 133-143.

\bibitem {03} L. Carlitz, {\it A note on Bernoulli and Euler polynomials of the second kind}, Scripta Math. 25 (1961), 323-330.

\bibitem{ 04}
R. Dere, Y. Simsek,
{\it Applications of umbral algebra to some special polynomials},
Adv. Stud. Contemp. Math. 22 (2012), no. 3, 433--438.


\bibitem {05} H. W. Gould, {\it Explicit formulas for Bernoulli numbers}, Amer. Math. Monthly 79(1972), 44-51.

\bibitem {06} G. Kim, B. Kim, J. Choi,  {\it The DC algorithm for computing sums of powers of consecutive integers and Bernoulli numbers}, Adv. Stud. Contemp. Math. 17 (2008), no. 2, 137-145.


\bibitem {07} D. S. Kim, N. Lee, J. N. Na, K. H. Park,  {\it Identities of symmetry for highter-order Euler polynomials in three variables (I)}, Adv. Stud. Contemp. Math. 22 (2012), no. 1, 51-74.

\bibitem {08} D. S. Kim, T. Kim, Y. H. Kim, D. V. Dolgy, {\it A note on Eulerian polynomials associated with Bernoulli and Euler numbers and polynomials}, Adv. Stud. Contemp. Math. 22 (2012), no. 3, 379-389.

\bibitem {09} D. S. Kim, T. Kim, S.-H Lee, Lee, S.-H Rim, {\it Some identities of Bernoulli, Euler and Abel polynomials arising from umbral calculus},  Adv. Difference Equ. 2013, 2013:15, 8 pp.

\bibitem {10} D. S. Kim, T. Kim, D. V. Dolgy, S.-H. Rim,  {\it Some new identities of Bernoulli, Euler and Hermite polynomials arising from umbral calculus}, Adv. Difference Equ. 2013, 2013:73.

\bibitem {11} D. S. Kim, T. Kim, S.-H Lee, {\it A note on poly-Bernoulli polynomials arising from umbral calculus}, Adv. Studies Theor. Phys. 7(2013), no. 15, 731-744.

\bibitem{7}
D. S. Kim, T. Kim, S. H. Lee,
 {\it Poly-Cauchy numvbers and polynomials with umbral calculus viewpoint},
Int. Journal of Math. Analysis, Vol. 7, 2013.

\bibitem {13} T. Komatsu, {\it On poly-Cauchy numbers and polynomials}, $http://carma.newcastle.edu.au$\\
$/alfcon/pdfs/Takao\_Komatsu-alfcon.pdf$


\bibitem {14} T. Komatsu, F. Luca,  {\it Some relationships between poly-Cauchy numbers and poly-Bernoulli numbers}, Annales Mathematicae et Informaticaw, 41(2013), 99-105.

\bibitem {15} T. Komatsu, {\it Poly-Cauchy numbers}, Kyushu J. Math., 67(2013), 143-153.

\bibitem {16} S. Roman, {\it The umbral calculus}, Pure and Applied Mathematics, 111, Academic Press, Inc. [Harcourt Brace Jovanovich, Publishers], New York, 1984. x+193 pp. ISBN: 0-12-594380-6.

\bibitem {17} S. Roman, G.-C. Rota,  {\it The umbral calculus}, Advances in Math. 27  (1978), no. 2, 95–188.

\end{thebibliography}
\end{document}